\documentclass[onecolumn]{autart}
\usepackage{graphicx,psfrag}          

\usepackage{amsmath}%
\usepackage{amsfonts}%
\usepackage{amssymb}%
\usepackage{array}%
\usepackage[hyphens]{url}

\graphicspath{{figures/}}

\newcommand{\range}{\operatorname{Im}}

\newcommand{\rank}{\operatorname{rank}}

\newcommand{\RR}{\ensuremath{\mathbb{R}}}

\newcommand{\N}{\mathcal{N}}

\newtheorem{assumption}{Assumption}
\newtheorem{remark}{Remark}

\def\varsetofrows{U}
\def\varperturbation{b}

\begin{document}

\begin{frontmatter}



\title{A Cyber Security Study of a SCADA Energy Management System: Stealthy Deception Attacks on the State Estimator\thanksref{footnoteinfo}}

\thanks[footnoteinfo]{This work was supported in part by the European Commission through the VIKING project, the Swedish Research Council, the Swedish Foundation for Strategic Research, and the Knut and Alice Wallenberg Foundation.}

\author[A]{Andr\'{e} Teixeira} 
\author[G]{Gy\"{o}rgy D\'{a}n} 
\author[A]{Henrik Sandberg}
\author[A]{Karl H. Johansson}

\address[A]{School of Electrical Engineering - Automatic Control, KTH~Royal~Institute~of~Technology}                                              
\address[G]{School of Electrical Engineering - Communication Networks, KTH~Royal~Institute~of~Technology}



\begin{abstract} 
The electrical power network is a critical infrastructure in today's society, so its safe
and reliable operation is of major concern. State estimators are commonly used
in power networks, for example, to detect faulty equipment and to optimally route power flows. The estimators are often located in control centers, to which large numbers of measurements are sent over unencrypted communication channels. Therefore cyber security for state estimators
becomes an important issue. In this paper we analyze the cyber security of state estimators in
supervisory control and data acquisition (SCADA) for energy management systems (EMS)
operating the power network. Current EMS state estimation algorithms have bad data detection
(BDD) schemes to detect outliers in the measurement data. Such schemes are based on
high measurement redundancy. Although these methods may detect a set of basic cyber
attacks, they may fail in the presence of an intelligent attacker. We explore the latter by
considering scenarios where stealthy deception attacks are performed by sending false information to the control center. We begin by presenting a recent framework that characterizes the attack as an optimization problem with the objective specified through a security metric and constraints corresponding to the attack cost. The framework is used to conduct realistic experiments on a state-of-the-art SCADA EMS software for a power network example with 14 substations, 27 buses, and 40 branches. The results indicate how state estimators for power networks can be made more resilient to cyber security attacks. 
\end{abstract}

\end{frontmatter}

\section{Introduction}
Examples of critical infrastructures in our society are the power, the gas and the water supply networks.
These infrastructures are operated by means of complex supervisory
control and data acquisition (SCADA) systems, which transmit
information through wide and local area networks to a control
center. Because of this fact, critical infrastructures are
vulnerable to cyber attacks,
see~\cite{kn:Johansson09,kn:Cardenas08b}. For a more recent
example that also received considerable media attention,
see~\cite{symantec10}.

\begin{figure}
  \begin{center}
   \includegraphics[width=0.5\hsize]{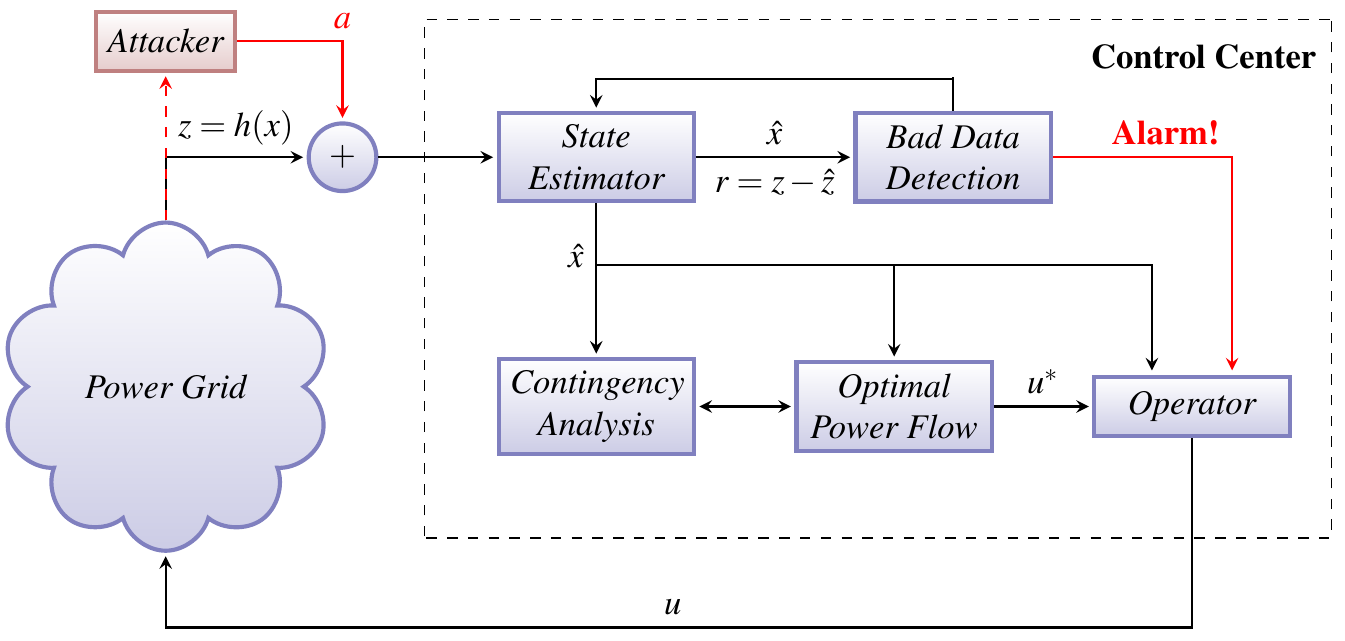}
  \end{center}
  \caption{The state estimator under a cyber attack\label{fig:SE_Attacked}}
\end{figure}

SCADA systems for power networks are complemented by a set of
application specific software, usually called energy management
systems (EMS). Modern EMS provide information support for a
variety of applications related to power network monitoring and
control. The power network state estimator (SE) is an on-line
application which uses redundant measurements and a network model
to provide the EMS with an accurate state estimate at all times.
The SE has become an integral tool for EMS, for instance for contingency analysis (CA) which, based on the state estimate, identifies the most severe consequences in case of hypothetical equipment outages.   
SCADA systems collect measurement
data from remote terminal units (RTUs) installed in various
substations, and relay aggregated measurements to the central
master station located at the control center. A simple schematic
picture of such a system is shown in Fig.~\ref{fig:SE_Attacked},
with measurements denoted by $z$.
 Several cyber attacks on SCADA systems operating power
networks have been reported, and major blackouts, such as the
August 2003 Northeast U.S. blackout, are due to the misuse of the
SCADA systems, see~\cite{kn:Blackout03}. As discussed
in~\cite{kn:Johansson09}, there are also several vulnerabilities
in the SCADA system architecture, including the direct tampering
of RTUs, communication links from RTUs to the control center, and
the IT software and databases in the control center.

Our work analyzes the cyber security of the SE in the SCADA system
of a power network. In current implementations of SE algorithms,
there are bad data detection (BDD) schemes
\cite{kn:Monticelli99,kn:Abur04} designed to detect random
outliers in the measurement data. Such schemes are based on high
measurement redundancy and are performed at the end of the state
estimation process. Although such methods may detect basic cyber
attacks on the measurements, they may fail in the presence of a
more intelligent attacker. It is well known that for so-called
\emph{multiple interacting bad data}, the BDD system can fail to
detect and locate the faulty measurements, see
\cite{kn:Monticelli99,kn:Abur04}. That an attacker can exploit
this fact has been pointed out in several recent papers, see
\cite{kn:Reiter09,kn:Sandberg10,kn:Bobba10}. For example, it has
been shown that an attacker with access to a model of the network
systematically can search for, and often find, simple undetectable
attacks. Returning to Fig.~\ref{fig:SE_Attacked}, this means it is
possible to compute data corruptions $a$ to measurements $z$ that
will not generate alarms in the control center. Such corruptions
are called \emph{stealthy deception attacks}.

In the work \cite{kn:Reiter09,kn:Sandberg10,kn:Bobba10}, it is
assumed that the attacker has a linear accurate model of the power
grid, and undetectability of the corruption $a$ is proven under
this assumption. The real power network is nonlinear, however, and a
nonlinear model is also typically implemented in the SE.
Therefore, it is not clear how a real SE will react to these
stealthy deception attacks. For example, how large can $a$ be
before the SE does no longer converge? In \cite{kn:Teixeira10}, we
have quantified how the SE residual can be bounded based on the
model error, but no tests on an actual system were performed. 

The main contribution of this paper is to test how sensitive a
state-of-the-art SCADA system SE is to stealthy deception
attacks. Maybe somewhat surprisingly, for the cases we have
studied, the attacks indeed pass undetected for very large
corruptions $a$. However, our analysis also shows that it is
possible to make these attacks much more difficult to perform by
allocating new sensors, or by securing some of them. Secure sensor
allocation has also been discussed in
\cite{kn:Bobba10,kn:Gyorgy10}.


The outline of the paper is as follows. In
Section~\ref{sec:Models} we present the theoretical concepts
behind state estimation in power networks. Results from previous
work are used in Section~\ref{sec:Framework} to develop the
analysis framework and some novel considerations regarding
limitations of linear attack policies are also given.
Section~\ref{sec:Results} contains the main contribution of this
paper, the description and results of practical experiments
conducted in a state-of-the-art SCADA/EMS software using the
previously mentioned framework. The conclusions are presented in Section~\ref{sec:conclusions}.

\section{Preliminaries}\label{sec:Models}
In this section we introduce the power network models and the theory behind the SE and BDD algorithms.

\subsection{Measurement model}
For an $N-$bus electric power network, the $n=2N-1$ dimensional state vector~$x$ is $(\theta^{\top},V^{\top})^{\top}$, where $V=(V_1,\ldots,V_{N})$ is the vector of bus voltage magnitudes and $\theta=(\theta_{2},\dots,\theta_{N})$ vector of phase angles. This state vector is the minimal information needed to characterize the operating point of the power network. Without loss of generality, we have considered bus~$1$ to be the reference bus, hence all phase-angles are taken relatively to this bus and $\theta_1 = 0$. The $m-$dimensional measurement vector $z$ can be grouped into two categories: (1) $z_{P}$, the active power flow measurements $P_{ij}$ from bus $i$ to $j$ and active power injection measurement $P_{i}$ at bus $i$, and (2) $z_{Q}$, the reactive power flow measurements $Q_{ij}$ from bus $i$ to $j$, reactive power injection measurement $Q_{i}$ and $V_{i}$ voltage magnitude measurement at bus $i$. The neighborhood set of bus $i$, which consists of all buses directly connected to this bus, is denoted by $N_i$. The power injections at bus $i$ are described by
\begin{equation*}
  \begin{array}{lcl}
    P_i & = & V_i\sum_{j\in N_i}V_j\left(G_{ij}\cos(\theta_{ij})+B_{ij}\sin(\theta_{ij})\right) \\
    Q_i & = & V_i\sum_{j\in N_i}V_j\left(G_{ij}\sin(\theta_{ij})-B_{ij}\cos(\theta_{ij})\right)
  \end{array} ,
\end{equation*}
and the power flows from bus $i$ to bus $j$ are described by 
\begin{equation*}
  \begin{array}{lcl}
    P_{ij} & = & V_i^2(g_{si}+g_{ij})-V_iV_j\left(g_{ij}\cos(\theta_{ij})+b_{ij}\sin(\theta_{ij})\right) \\
    Q_{ij} & = & -V_i^2(b_{si}+b_{ij})-V_iV_j\left(g_{ij}\sin(\theta_{ij})-b_{ij}\cos(\theta_{ij})\right)
  \end{array} ,
\end{equation*}
where $\theta_{ij}=\theta_i-\theta_j$ is the phase angle difference between bus $i$ and $j$, $g_{si}$ and $b_{si}$ are the shunt conductance and susceptance of bus $i$, $g_{ij}$ and $b_{ij}$ are the conductance and susceptance of the branch from bus $i$ to $j$, and $Y_{ij}=G_{ij}+jB_{ij}$ is the $ij$th entry of the nodal admittance matrix. 
More detailed formulas relating measurements~$z$ and state~$x$ may be found in~\cite{kn:Abur04}.

Assuming that the model parameters and the network topology are exact, the nonlinear measurement model for state estimation is defined by
\begin{align}\label{eq:measmodel}
z=h(x)+\epsilon,
\end{align}
where $h(\cdot)$ is the $m-$dimensional nonlinear measurement function that relates measurements to states and is assumed to be twice continuously differentiable, $\epsilon=(\epsilon_{1},\ldots,\epsilon_{m})^{\top}$ the zero mean measurement error vector, and usually $m\gg n$ meaning that there is high measurement redundancy. Here $\epsilon_{i}$ are independent Gaussian variables with respective variances~$\sigma_{i}^2$ indicating the relative uncertainty about the $i-$th measurement and thus we have $\epsilon\sim\mathcal{N}(0,R)$ where $R=\text{diag}(\sigma_{1}^2,\ldots,\sigma_{m}^2)$ is the covariance matrix.

\subsection{State Estimator} 
The basic SE problem is to find the best $n$-dimensional state $x$ for the measurement model \eqref{eq:measmodel}
in a weighted least square (WLS) sense. Defining the residual vector $r(x)=z-h(x)$, we can write the WLS problem as
\begin{equation}\label{eq:SE_eq}
\begin{aligned}
			 \min_{x\in\mathbb{R}^n} J(x)&=\frac{1}{2}r(x)^{\top}R^{-1}r(x) \\
\text{such that }	 g(x)&=0\\
			 s(x)&\leq0,
\end{aligned}
\end{equation}
where the inequality constraints generally model saturation limits, while the equality constraints are used to include target setpoints and to ensure physical laws such as zero power injection transition buses, e.g., transformers, and zero power flow in disconnected branches. Thus data used in the equality constraints is often seen as \emph{pseudo-measurements}. For sake of simplicity, we will present the solution to the unconstrained optimization problem. More detailed information on the solution of~\eqref{eq:SE_eq} may be found in~\cite{kn:Abur04} and~\cite{kn:Monticelli99}.

The unconstrained WLS problem is posed as
\[\min_{x\in\mathbb{R}^n} J(x)=\frac{1}{2}r(x)^{\top}R^{-1}r(x).\]
The SE yields a \emph{state estimate} $\hat x$ as a minimizer to this problem. 
The solution $\hat x$ can be found using the \emph{Gauss-Newton} method which solves the so called \emph{normal equations}:
\begin{align}\label{eq:normalequations}
\left(H^{\top}(x^k)R^{-1}H(x^k)\right) (\Delta x^k) = H^{\top}(x^k)R^{-1}r(x^k),  
\end{align}
for $k=0,1,\ldots$, where \[H(x^k):=\frac{d h(x)}{d x}|_{x=x^k}\] is called the Jacobian matrix of the measurement model $h(x)$. 
For an observable power network, the measurement Jacobian matrix $H(x^k)$ is full column rank.
Consequently, the matrix $\left(H^{\top}(x^k)R^{-1}H(x^k)\right) $ in~\eqref{eq:normalequations} is positive definite and the Gauss-Newton step generates a descent direction, i.e., for the direction $\Delta x^k=x^{k+1}-x^k$ the condition $\nabla J(x^k)^{\top}\Delta x^k<0$ is satisfied.

\begin{remark}
 Henceforth we consider the covariance matrix $R$ to be the identity matrix, i.e., all measurements have unitary weights. The framework and results presented in the next sections can be easily extended to the more general case, see~\cite{kn:Teixeira10}.
\end{remark}

For notational convenience, throughout the next sections we will use $H(x^k)$ as $H$, $\Delta x^k$ as $\Delta x$, and $r(x^k)=z-h(x^k)$ as $r$.


\subsection{Decoupled State Estimation}

A useful observation in electric power networks is that of active-reactive decoupling, i.e., the active measurements $z_{P}$ (resp. reactive measurement $z_{Q}$) predominantly affect the phase angles $\theta$ (resp. the voltage magnitudes $V$). In the decoupled state estimation, the approximate values of the corrections $\Delta \theta$ and $\Delta V$ are then not computed simultaneously, but independently~\cite{kn:Wu90_Survey}.

Following~\eqref{eq:normalequations}, the correction to state estimate $\Delta x=(\Delta \theta^{\top}, \Delta V^{\top})^{\top}$ at each iteration can be obtained from the weighted measurement residual $r = (r_{P}^\top, r_{Q}^\top)^\top$ as the solution to the overdetermined system
\begin{align}\label{eq:matnormal}
\begin{pmatrix}
H_{P\theta}  & H_{PV} \\  H_{Q\theta}  & H_{QV}
\end{pmatrix} \begin{pmatrix}
\Delta \theta \\ \Delta V
\end{pmatrix}=\begin{pmatrix} r_{P} \\ r_{Q} \end{pmatrix},
\end{align}
where the submatrices $H_{P\theta}$ and $H_{PV}$ correspond to active measurements and $H_{Q\theta}$ and $H_{QV}$ correspond to reactive measurements. The traditional version of fast decoupled state estimation is based on the following decoupled normal equations, where the coupling submatrices $H_{PV}$ and $H_{Q\theta}$ have been set to zero: 
\begin{align}\label{eq:fastdecoupled}
\begin{split}
\Delta \theta^{k}=H_{P\theta}^{\dagger} r_{P}(\theta^k, V^k),\\
\Delta  V^{k}= H_{QV}^{\dagger}r_{Q}(\theta^k, V^k).
\end{split}
\end{align}
Equations~\eqref{eq:fastdecoupled} are alternately solved for $\Delta \theta^k$ and $\Delta V^k$, where the mismatches $r_{P}$ and $r_{Q}$ are evaluated at the latest estimates. The submatrices $H_{P\theta}$ and $H_{QV}$ are evaluated at \emph{flat start} and branch series resistances are ignored in forming $H_{P\theta}$. By flat start we mean the power network's state in which all voltage magnitudes are 1 pu and all phase angles are 0.

\subsection{Bad Data Detection}
\label{sec:BDD}

The measurement residual when random bad data is present is characterized as follows.
Assume there are no measurement errors, i.e. $z=h(x)$, and that the SE has converged through the Gauss-Newton method. Recalling that $r(\hat{x})=z-h(\hat{x})$, from~\eqref{eq:normalequations} we see that the estimate sensitivity matrix is given by $\frac{\partial \hat{x}}{\partial z} = (H^{\top}H)^{-1}H^\top$. Furthermore, we conclude that the weighted residual sensitivity matrix is $\frac{\partial r}{\partial z}=I-\frac{\partial h(\hat{x})}{\partial \hat{x}} \frac{\partial \hat{x}}{\partial z} =I-H(H^{\top}H)^{-1}H^\top$. Thus for small measurement errors $\epsilon\sim\N(0,I)$ we have the following weighted measurement residual
\begin{equation}
 \label{eq:weighted_res}
 r=S\epsilon,
\end{equation}
where $ S=I-H(H^{\top}H)^{-1}H^\top$.

Through BDD the SE detects measurements corrupted by errors whose statistical properties exceed the presumed standard deviation or mean. This is achieved by hypothesis tests using the statistical properties of the weighted measurement residual~\eqref{eq:weighted_res}. We now introduce one of the BDD hypothesis tests widely used in practice, the \emph{largest normalized residual test}.

\subsubsection{Largest normalized residual test}
From~\eqref{eq:weighted_res}, we note that $r\sim\N(0,\Omega)$ with $\Omega = S$. Now consider the normalized residual vector
\begin{equation}\label{eq:rNdefn}
 r^N=D^{-1/2}r,
\end{equation}
with $D\in\mathbb{R}^{m\times m}$ being a diagonal matrix defined as $D=\operatorname{diag}(\Omega)$. In the absence of bad data each element $r^N_i,\ i=1,\dots,m$ of the normalized residual vector then follows a normal distribution with zero mean and unit variance, $r^N_i\sim\N(0,1),\ \forall i=1,\dots,m$. Thus, bad data could be detected by checking if $r^N_i$ follows $\N(0,1)$. This can be posed as a hypothesis test for each element $r^N_i$
\begin{align*}
 H_0: \mathbb{E}\left\{r^N_i\right\} = 0, \quad H_1: \mathbb{E}\left\{|r^N_i|)\right\} > 0.
\end{align*}
For this particular case, as shown in~\cite{kn:Monticelli99}, the largest normalized residual (LNR) test corresponds to a threshold test where the threshold $\tau$ is computed for a given false alarm rate and $H_0$ is accepted if
\begin{equation}\label{eq:LNR}
 \|D^{-1/2}r\|_{\infty} \leq \tau,
\end{equation}
and rejected otherwise.

\section{Stealthy deception attacks}\label{sec:Framework}
Using the theory and models described in the previous section, we present the framework used throughout the next sections to study the cyber security of SCADA EMS software and algorithms.

\subsection{Attacker Model}
\label{sec:AttackerModel}
The goal of a stealthy deception attacker is to compromise the telemetered measurements available to the SE such that: 1) The SE algorithm converges; 2) The attack remains undetected by the BDD scheme; and 3) For the targeted set of measurements, the estimated values at convergence are close to the compromised ones introduced by the attacker.


Let the corrupted measurement be denoted $z^a$. We assume the following additive attack model
\begin{align}\label{eq:attackmodel}
z^a=z+a,
\end{align}
where $a\in\RR^m$ is the attack vector introduced by the attacker, see also Fig.~\ref{fig:SE_Attacked}. The vector~$a$ has zero entries for uncompromised measurements. Under attack, the normal equations~\eqref{eq:normalequations} give the estimates
\begin{align}\label{eq:normal_eq_attacked}
 \tilde x^{k+1} = \tilde x^{k}+\left(H^{\top}(\tilde x^k)H(\tilde x^k)\right)^{-1}H^{\top}(\tilde x^k)r^a(\tilde x^k),
\end{align}
for $k=0,1,\dots$, where $\tilde x^k$ is the \emph{biased} estimate at iterate $k$, and $r^a(\tilde x^k):=z^a-h(\tilde x^k)$. If the local convergence conditions hold, then these iterations converge to $\hat{x}^a$, which is the biased state estimate resulting from the use of $z^a$. Thus, the convergence behavior can be expressed as the following statement:
\begin{enumerate}\label{en:obj1}
\item[1)] The sequence $\{\tilde x^0,\tilde x^1,\dots\}$ generated by~\eqref{eq:normal_eq_attacked} converges to a fixed point~$\hat x^{a}$.
\end{enumerate}
We will occasionally use the notation $\hat x^a(z^a)$ to emphasize the dependence on $z^a$. 

The BDD scheme for SE is based on a threshold test. Thus the attacker's action will be undetected by the BDD scheme provided that the following condition holds:
\begin{enumerate}\label{en:obj2}
\item[2)] The measurement residual under attack $r^a:=r(\hat x^a)=z^a-h(\hat x^a)$, satisfies the condition~\eqref{eq:LNR}.
\end{enumerate}
  
Finally, consider that the attacker aims at corrupting measurement $i$. This means the attacker would like the estimated measurement $\hat z^a_i:=h_i(\hat x^a(z^a))$ to be equal to the actual corrupted measurement $z^a_i$. Therefore, we arrive at the following condition which will additionally govern the synthesis of attack vector~$a$:

\begin{enumerate}\label{en:obj3}
\item[3)] The attack vector $a$ is chosen such that $|z^a_i-\hat z^a_i|=0$.
\end{enumerate}

The aim of a stealthy deception attacker is then to find and apply an attack $a$ that satisfies conditions 1), 2), and 3). This problem can be posed as
  \begin{equation}\label{eq:valid_attack}
  \begin{aligned}
& \text{find } a \\
\text{s.t. } & a\in\mathcal{G}\cap\mathcal{C} \cap\mathcal{U}\ ,
   \end{aligned}
  \end{equation}
where $\mathcal{G}$ is the set of goals in condition 3), $\mathcal{C}$ the set of constraints ensuring condition 1) is met and that no protected or pseudo-measurements are corrupted, and $\mathcal{U}$ the set of stealthy attacks satisfying condition 2). 


\subsection{Security Metric}\label{sec:security_metric}
 
 In general a stealthy attack requires the corruption of more measurements than the targeted one, see~\cite{kn:Reiter09} and \cite{kn:Sandberg10}. Such requirement relates to the fact that a stealthy attack must have the attack vector $a$ fitting the measurement model. 
 
 Considering that the system's state is $x^*$ and the attacks are sufficiently small, the measurement model can be linearized around $x^*$, obtaining
  \begin{equation}\label{eq:linear_model}
   z=\frac{d h(x)}{d x}|_{x=x^*}(x^* + c)=H(x^*+c),
  \end{equation}
 where $c$ is the perturbation added to $x^*$. Previous results show that the class of stealthy attacks for this linear model is characterized by $a\in\range(H)$, which is equivalent to have $a=Hc$, for some $c\neq0$. Based on this linear model, we present a security metric $\alpha_k$ for each measurement $k$. This metric corresponds to the minimum cost of a valid attack satisfying~\eqref{eq:valid_attack} and targeting to corrupt measurement $k$ by adding it one unit, i.e., $a_k=1$. It is computed by solving the problem
    \begin{equation}
  \begin{aligned}
\alpha_k=& \min_a \|a\|_0 \\
\text{s.t. } & a\in\mathcal{G}\cap\mathcal{C}\cap\mathcal{U}\ ,
   \end{aligned}
  \end{equation}
where here $a\in\mathcal{G}$ corresponds to having $a_k=1$, $a\in\mathcal{C}$ to $a_i=0\,, \forall i$ if measurement $i$ is a pseudo-measurement, and $a\in\mathcal{U}$ to $a\in\range(H)$. Note that $\|\cdot\|_0$ is a pseudo-norm corresponding to the cardinality, i.e., number of non-zero entries, of the argument. Hence the cost $\alpha_k$ corresponds to the minimum cardinality of a valid attack $a$, i.e., the minimum number of sensors needed to be corrupted.

This metric can also be extended to cases where by compromising a single measurement, the attacker gains access to other measurements without additional cost. For a more detailed discussion on this metric and efficient algorithms to compute it, see~\cite{kn:Gyorgy10}.

\subsection{Limitations of Linear Policies}\label{sec:limitations}

In this section we comment on the limitations of the linear attack policies described in Section~\ref{sec:security_metric}. Recall that the core of the linear policies is to have $a\in\range(H)$. Two main limitations arise from this policy, one related to the fact that $H$ is obtained for given operating conditions and the other related to saturation limits not considered by the linear measurement model~\eqref{eq:linear_model}. We briefly discuss about these limitations.

\subsubsection{Varying operating conditions}

The power network is a dynamical system and its state is frequently changing. Thus it might be the case that the attacker has previously obtained a linear model $\tilde{H}$ for a state $\tilde{x}$ and the attack is performed only when the system is in a different state $x^*$ where the linear approximation is $H$. 
Hence for small attack vectors or for cases where $x^*\approx \tilde{x}$, the residual will be small and the attack may pass undetected. For larger attacks, however, this might not hold. These scenarios can be analyzed using the framework presented in~\cite{kn:Teixeira10}, where it is considered that the attacker has an inaccurate model $\tilde H$.

One interesting fact to observe is that, under certain assumptions, the attack vector $a$ and the security metric $\alpha$ will be the same, independently of the system's state. We now present a useful lemma and the required assumptions and formulate this result.

\begin{lem}\label{lemma:rank}
Consider an optimal attack $a^*$ that is undetectable with respect to $H$, 
i.e., $a^*\in \range(H)$ and has minimum cardinality. Denote by $\varsetofrows$ the set of measurements not
affected by $a^*$, i.e., $a^*_i=0\ \forall i\in\varsetofrows$ and
$a^*_i\not=0\ \forall i\not\in\varsetofrows$. Let the Jacobian matrix be partitioned as $H=[H_{\varsetofrows}^\top\  H_{\bar{\varsetofrows}}^\top]^\top$.
Recalling that $n$ is the number of states and that $\rank(H)=n$, then $\rank(H_\varsetofrows)=n-1$, and for every $i\not\in\varsetofrows$
we have $\rank(H_{\varsetofrows\cup\{i\}})=n$.
\end{lem}
\begin{pf}
 See~\cite{kn:Gyorgy10}.{\ \rule{0.5em}{0.5em}}
\end{pf}

\begin{assumption}\label{assump:nonzero}
 For any measurement element $z_i$ we have $\frac{\partial h_i(x^*)}{\partial x_j}=0$ if and only if $\frac{\partial h_i(\tilde{x})}{\partial x_j}=0$, for all $j=1,\dots,n$.
\end{assumption}
\begin{prop}
Denote $\mathcal{F}$ as the set of power flow measurements. Consider $H=H_{P\theta}(\theta)|_{\theta=\theta^*}$, $\tilde{H}=H_{P\theta}(\theta)|_{\theta=\tilde{\theta}}$, and let $a$ be a stealthy deception attack vector with respect to $H$. Denote the set of measurements not corrupted by $a$ as $\varsetofrows$. Then for all the line parameter perturbations and state changes from $\theta^*$ to $\tilde{\theta}$ that satisfy Assumption~\ref{assump:nonzero} and do not affect measurements $j\in\mathcal{F}\cap\bar{\varsetofrows}$, we have that $a$ is also a stealthy attack with respect to $\tilde{H}$.
\end{prop}

\begin{pf}
In the following we consider the matrix $H_{P\Theta}$
and use $H$ instead to simplify the notation.
We consider a perturbation in the linear model, 
e.g., due to varying operating conditions, such
that for a measurement $k\in\mathcal{F}$ corresponding to a 
transmission line we have $\tilde{H}_k=\varperturbation H_k$.
Let us denote the buses at the two ends of the transmission
line by $k_1$ and $k_2$. For the power injection at
bus $k_1$ (closest to measurement $k$) we have
$\tilde{H}_{k_1}=H_{k_1}+(\varperturbation -1)H_{k}$,
for bus $k_2$ we have 
$\tilde{H}_{k_2}=H_{k_2}-(\varperturbation -1)H_{k}$.
In the following we show that if $k\in\mathcal{F}\cap\varsetofrows$, then for every $a\in\range(H)$ there is
$\tilde{a}\in\range(\tilde{H})$ such that $a=\tilde{a}$. Since this also holds for minimum cardinality attack
vectors, we have that the security metric $\alpha$ is the same for both linearized models.

For the case when $k\in\varsetofrows$ we can prove the proposition
by performing elementary row operations on $\tilde{H}_{\varsetofrows}$.
If $k_1\in\varsetofrows$ we subtract $(\varperturbation-1)/\varperturbation\tilde{H}_{k}$
from $\tilde{H}_{k_1}$. We proceed similarly for $k_2$.
Finally, we divide $\tilde{H}_{k}$ by $\varperturbation$.
Clearly, after these operations we obtained $H_{\varsetofrows}$,
which, following Lemma~\ref{lemma:rank}, proves that $\rank(\tilde{H}_{\varsetofrows})=n-1$.
Observe that since we used elementary row operations, the kernel of
$\tilde{H}_{\varsetofrows}$ is the same as that of $H_{\varsetofrows}$.
Consequently, the same attack vectors can be used despite
the perturbation of the model, i.e., $a\in\range(H)\cap\range(\tilde{H})$. \ \rule{0.5em}{0.5em}
\end{pf}

Hence we conclude that if the state or parameter perturbations do not affect power flow measurements compromised by the attack before the change, then the same attack vector is still valid after the parameter or state change. Note that this is ensured if the measurements affected by the parameter change are far from the region of the network where the attack is performed. Thus this indicates that attacks can be performed locally in the power network. Additionally, in this case the security metric is the same for both linearized models and there is no need to recompute it again for the different state. 

\subsubsection{Saturation limits}

The linear measurement model~\eqref{eq:linear_model} is obtained by linearizing the nonlinear model~\eqref{eq:measmodel} at a given state $x^*$ and so the linear model only approximates well the nonlinear one in a region close to $x^*$. Furthermore, the more nonlinear the function is around $x^*$, the smaller is the region where the approximation is valid. This fact is particularly important when the saturation occurs. From considering~\eqref{eq:linear_model} alone we do not have any limits on the size of $a=Hc$. However, the nonlinear model clearly shows that the measurements have saturation limits. For instance, disregarding the line and shunt conductances $g_{ij}$ and $g_{si}$ we have $P_{ij}=-V_iV_jb_{ij}\sin(\theta_{ij})$, where we see that the theoretical maximum for this power flow is given by $P_{ij}=-V_iV_jb_{ij}$. Hence for a stealthy attack it is not enough to require that $a\in\range(H)$ but it is also essential to impose saturation limits on the attacked measurements. These limits could be formulated as inequality constraints and included in the set of constraints $\mathcal{C}$, in general reducing the set of valid attacks.

\section{Experiments on the SCADA EMS system}\label{sec:Results}
During the previous sections we have mentioned recent work where the authors analyzed stealthy deception attacks on SE based on linearized models.
However, the results obtained so far do not clarify how sensitive the real SCADA EMS software is to these attacks or if a system operator should even care about these scenarios. In this section we present the results obtained by carrying out a stealthy deception attack on a real SCADA EMS software. We hope, by analyzing these results, to answer the previous open questions and also to provide recommendations to increase the security of SCADA EMS software against deception attacks. Before analyzing the results, we briefly describe the experimental setup.

\subsection{Experimental Setup}
\begin{figure*}
 \centering
 \includegraphics[width=0.6\hsize]{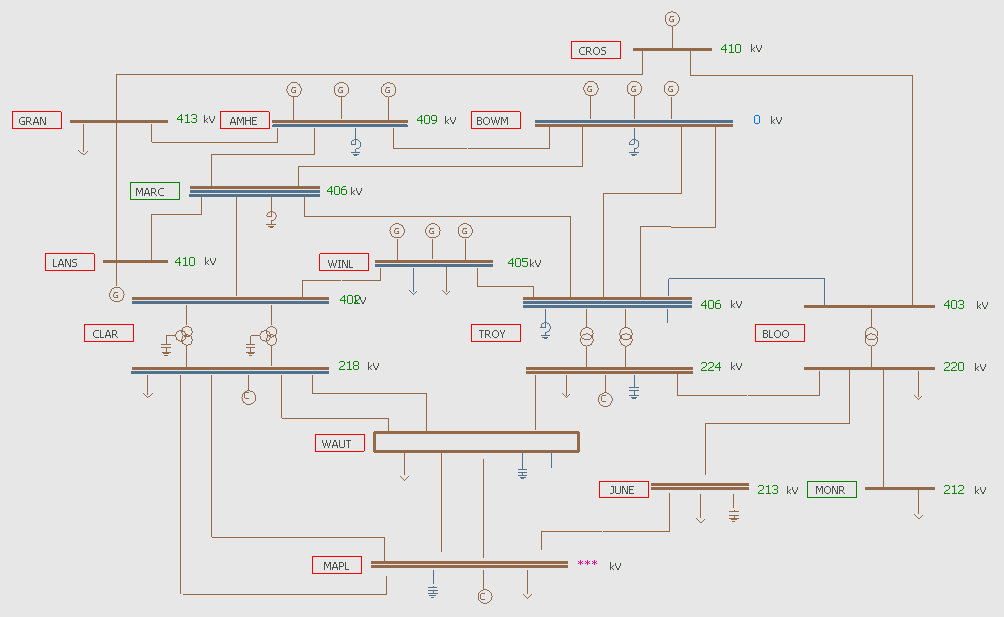}
\caption{Power network considered in the experiment, where the gray and brown colors represent unenergized and energized objects, respectively. \label{fig:PowerNetwork}}
\end{figure*} 

The software was supplied with the virtual network presented in Fig.~\ref{fig:PowerNetwork}, similar to the IEEE 39 bus network. The power network in Fig.~\ref{fig:PowerNetwork} consists of 14 substations and the bus-branch model has 27 buses and 40 branches. Several measurements are available at each substation, such as voltage magnitudes, active and reactive power flows and injections, and transformer tap change positions. This data is kept in the software database. We have used a console-based static network simulator to carry out the data corruption by directly changing the measurement data in the database.
The presented results thus relate to data corruption attacks and the consequences of such attacks on the EMS software components.

Specific EMS components, such as SE and BDD, are configured with unitary weights for all the measurements. The SE solves the nonlinear weighted least-squares problem using the fast-decoupled algorithm with equality constraints, while the BDD algorithm uses the LNR test. Both approaches correspond to standard algorithms presented in Section~\ref{sec:Models}.

As described in previous sections, some information about the power network is needed to compute stealthy deception attacks. Here we consider a particular class of such information, namely the bus-branch model 
of the network. In this experiment, we exported this information to MATLAB using the MATPOWER toolbox,~\cite{kn:MATPOWER}. A simplified attack was considered in which only the DC model of the network was used. This corresponds to including only active power measurements in the set of corrupted measurement data, disregarding the reactive measurements, not taking into account the current operating state of the system, not the coupling between the active power and the voltage magnitudes, and not the line conductances or the shunt admittances. Hence all voltage magnitudes were assumed to be 1 pu and the phase-angles 0. Only a simplified version of the $H_{P\theta}$ submatrix in~\eqref{eq:matnormal} was used, hereby denoted $H_{DC}$.


The algorithm in~\cite{kn:Gyorgy10} was used to compute the security metrics for each measurement. Information regarding which measurements were assumed to be tamper-proof was taken in account. Such measurements correspond to pseudo-measurements, which are considered as equality constraints in the optimization problem, and are often based on physical principles, see~\cite{kn:Monticelli99}. 

\begin{figure}[ht]
 \centering
 \includegraphics[viewport=81   320   529   471, width=0.7\hsize]{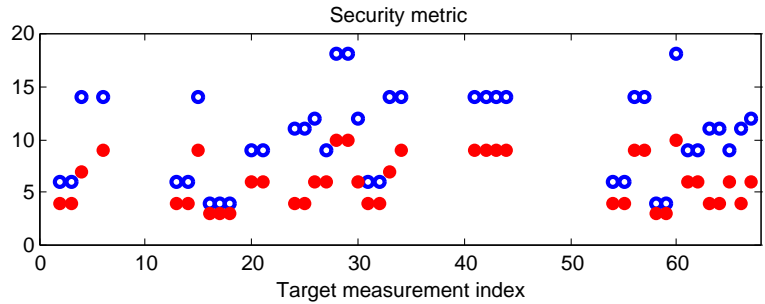}
 \caption{Security metrics for each measurement $k$: $\alpha_k$ (red full circles) was computed taking in account which measurements are available in the SCADA EMS, while $\bar{\alpha}_k$ (blue rings) was computed assuming that all possible measurements are being taken. Both represent the minimum number of measurements needed to perform a stealthy attack on the target measurement $k$.\label{fig:security}}
\end{figure}

The result is presented in Fig.~\ref{fig:security}. Given the current configuration of the SCADA EMS, specifically which measurements are available, we computed the security metric $\alpha_k$ (the red full circles) as defined in Section~\ref{sec:security_metric}. We see that the result is very heterogeneous, since around a third of the measurements has low values between 3 and 4, another third has values between 6 and 10, and others, which are not depicted, have values greater than 20 or even $\alpha_k=\infty$. Recalling that $\alpha_k$ is the minimum number of measurements needed to perform a stealthy attack on measurement $k$, we conclude that measurements with low $\alpha_k$ are easily attacked while the ones with $\alpha_k=\infty$ are fully protected.

Increasing the redundancy of the system by adding more measurements to the SCADA system increases the security level, as we see by looking at how $\bar\alpha_k$ is larger than $\alpha_k$. However, note that this does not guarantee full protection, as all measurements with finite $\alpha_k$ still have finite $\bar\alpha_k$. 
\subsection{Attack Scenario}

\begin{table}
 \caption{Example: adding 100MW to target measurement $33$\label{tab:attack}}
 \centering
 \begin{tabular}{cccc} \hline
 Measurement	& Normalized	& Correct value	& False value	\\
 index,	$k$	& attack, $\bar{a}_k$	&  (MW), $z^*_k$	& (MW), $z^a_k$	\\ \hline
 4 	& -1 		& 1005.7041 	& 905.7042\\
 21	& -0.7774	& 157.8541 	& 80.1103\\
 24	& 0.9665	& 507.7171 	& 604.3638\\
 27	& 2.7439	& 40.0006		& 314.3911\\
 \textbf{33} & \textbf{1}	& \textbf{-14.7971} & \textbf{85.2029}\\
 62	& 0.7774	& -123.3764	& -45.6327\\
 104	& -0.9665	& -334.8826	& -431.5293
 \end{tabular}
\end{table}
To conduct our experiment we considered measurement number 33, corresponding to the active power flow on the tie-line between TROY and BLOO substations, to be the target measurement. This means that the attacker's goal is to change this power flow measurement value as he/she wishes. In order to do so without being detected, the attacker needs to perform a coordinated attack in which he/she corrupts the value of other power measurements. Following the framework presented in Section~\ref{sec:Framework}, the set of such malicious changes is encoded in the attack vector $a$, which is then added to the true measurement vector $z$. The corrupted measurement vector $z^a=z+a$ is the one used by the SE.

Using $H_{DC}$, we computed the additive normalized attack vector required to stealthily change the target measurement by 1 MW, presented in Table~\ref{tab:attack}. As seen in Fig.~\ref{fig:security}, such attack only corrupts 7 measurements in total, which are taken from 5 substations, namely TROY, BLOO, JUNE, MONR, and CROS, all situated in the right side of Fig.~\ref{fig:PowerNetwork}. Hence we see that to stealthily attack a single measurement, a local coordinated attack suffices, even for such a large system. Additionally, as discussed in~\cite{kn:Gyorgy10}, note that usually all measurements within a given substation are gathered at a single RTU. This means that by breaking into the substation's RTU the attacker gains access to all those measurements, so we can argue that although 7 measurements need to be corrupted, only 5 RTU's need to be compromised.

\subsection{Experimental Results}
The normalized attack vector $\bar{a}$, whose non-zero entries are shown in Table~\ref{tab:attack}, was used to corrupt the measurement data according to the attacker's objective. For instance, in Table~\ref{tab:attack} we can see the correct value of the compromised measurements, denoted by $z^*$, and the false values sent to the control center, $z^a$, when the objective was to induce a bias of 100MW in the target measurement, having $z^a=z^*+100\bar{a}$.

\begin{figure}
 \centering
 \includegraphics[viewport=81   227   529   564,width=0.5\hsize]{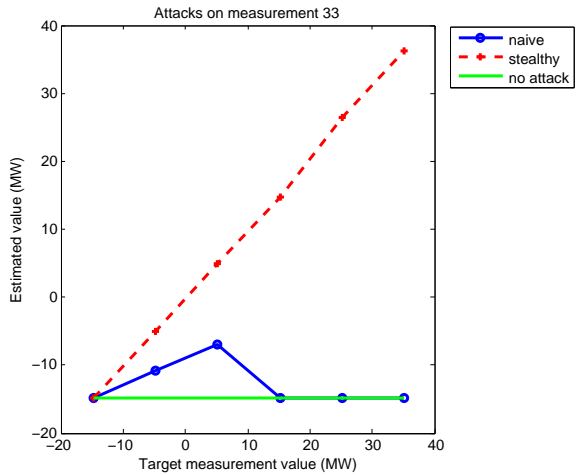}
 \caption{Stealthy deception attack\label{fig:small_attack}}
\end{figure}

In Fig.~\ref{fig:small_attack} we show the results obtained by performing stealthy deception attacks as described before and naive deception attacks where only the target measurement is compromised. In both cases, the bias in the target measurement was sequentially increased by 10MW at each step. From these results we see that the naive attack was undetected up to a bias of 20MW, while for bias above 30MW this attack was detected and the compromised measurement removed. The coordinated stealthy attack, however, remained undetected for all the bias values showed in the figure. Furthermore we see that the naive attack did not influence the estimate as much as the stealthy one, for which the relationship between the false and the estimated values is an almost unitary slope.

Table~\ref{tab:attack_data} shows the results obtained for large bias, where the attacks were performed sequentially with steps of 50MW. We observe that the stealthy attacks were successful with no BDD alarm triggered up to a bias of 150MW, beyond which the SE no longer converged.

\begin{table}
 \caption{Results from the stealthy attack for large bias\label{tab:attack_data}}
 \centering
 \begin{tabular}{ccccc} \hline
  Target bias, & False value & Estimate & \#BDD & \#CA\\
  $a_{33}$&  (MW), $z^a_{33}$ &(MW), $\hat{z}^a_{33}$ & Alarms & Alarms\\ \hline
 0& -14.8& -14.8& 0 & 2\\
 50&35.2& 36.2& 0 & 2\\
 100&85.2& 86.7& 0 & 10\\
 150&135.2& 137.5& 0 & 27\\
 200&185.2& - & - & -\\
 \end{tabular}
\end{table}


Although the SE did not converge for attacks above 200MW, it is still surprising to see that attacks based on the linearized model as large as 150MW are successful. To better understand what such quantity indicates, note that the nominal value of the targeted tie-line is 260MW. Thus the attack was able to induce a bias of more than 50\% of the nominal value, which reveals that the SCADA EMS software is indeed sensitive to stealthy deception attacks. Furthermore, notice that the number of warnings given by the CA component increase with the size of the attack. Whether or not this trend is related to the fact that the coupling terms have been neglected is still a matter for future analysis. Nevertheless, note that the increased number of CA warnings could lead the operator to take corrective actions. Therefore, we conclude that operators and utilities should care about these scenarios.

We also want to highlight that these results were achieved with a simplified linear model where several parameters, including the correct operating conditions and cross-coupling effects between active and reactive measurements, were disregarded. However in these scenarios we assumed the attacker had a large amount of resources such as a rather detailed knowledge regarding the network model, the available measurements, and the pseudo-measurements, and access to several RTUs. Most likely, an attacker with such resources could perform more devastating attacks on the power network than the ones considered here.

\section{Conclusions}\label{sec:conclusions}
In this paper we presented a comprehensive framework to analyze and study a class of stealthy deception attacks specifically targeting the SE component of SCADA EMS software through measurement data corruption. This framework provides attacker and attack cost models, possible attack synthesis policies, and system security metrics. The system security metric can be used by the utility to strengthen the security of the system by allocation of new sensors. Some limitations of the linear attack policies were briefly discussed. To validate this framework, we conducted a set of deception attacks to a state-of-the-art SCADA EMS software. The results obtained by this experiment show that computations based on linear models of the system provide valid attacks that successfully corrupt the target measurements without triggering any BDD alarms. The results also indicate that linear models can be used for large attacks as well, although otherwise expected. Additionally, we showed that besides the measurement model, information concerning pseudo-measurements and saturation limits is needed for a successful stealthy attack. This study also shows that improved BDD schemes and methods to ensure measurement and data protection are desirable.

\begin{ack}                               
The authors would like to thank
Mr.~Moustafa Chenine and Mr.~Nicholas Honeth
for their helpful technical support.
\end{ack}

\bibliographystyle{ieeetr}        
\bibliography{references_pss}           

\begin{thebibliography}{10}

\bibitem{kn:Johansson09}
A.~Giani, S.~Sastry, K.~H. Johansson, and H.~Sandberg, ``The {VIKING} project:
  an initiative on resilient control of power networks,'' in {\em Proc. 2nd
  Int. Symp. on Resilient Control Systems}, (Idaho Falls, ID, USA), pp.~31--35,
  Aug. 2009.

\bibitem{kn:Cardenas08b}
A.~C{\'a}rdenas, S.~Amin, and S.~Sastry, ``Research challenges for the security
  of control systems.,'' in {\em Proc. 3rd USENIX Workshop on Hot topics in
  security}, July 2008.

\bibitem{symantec10}
Symantech, ``Stuxnet introduces the first known rootkit for industrial control
  systems,'' August 6th 2010.

\bibitem{kn:Blackout03}
{U.S.-Canada Power System Outage Task Force}, ``Final report on the {A}ugust
  14th blackout in the {U}nited {S}tates and {C}anada,'' tech. rep., April
  2004.

\bibitem{kn:Monticelli99}
A.~Monticelli, {\em State Estimation in Electric Power Systems: A Generalized
  Approach}.
\newblock Kluwer Academic Publishers, 1999.

\bibitem{kn:Abur04}
A.~Abur and A.~Exposito, {\em Power System State Estimation: Theory and
  Implementation}.
\newblock Marcel-Dekker, 2004.

\bibitem{kn:Reiter09}
Y.~Liu, M.~K. Reiter, and P.~Ning, ``False data injection attacks against state
  estimation in electric power grids,'' in {\em Proc. 16th ACM Conf. on
  Computer and Communications Security}, (New York, NY, USA), pp.~21--32, 2009.

\bibitem{kn:Sandberg10}
H.~Sandberg, A.~Teixeira, and K.~H. Johansson, ``On security indices for state
  estimators in power networks,'' in {\em Preprints of the First Workshop on
  Secure Control Systems, CPSWEEK 2010}, Apr. 2010.

\bibitem{kn:Bobba10}
R.~Bobba, K.~M. Rogers, Q.~Wang, H.~Khurana, K.~Nahrstedt, and T.~Overbye,
  ``Detecting false data injection attacks on {DC} state estimation,'' in {\em
  Preprints of the First Workshop on Secure Control Systems, CPSWEEK 2010},
  Apr. 2010.

\bibitem{kn:Teixeira10}
A.~Teixeira, S.~Amin, H.~Sandberg, K.~H. Johansson, and S.~S. Sastry, ``Cyber
  security analysis of state estimators in electric power systems,'' in {\em
  Proc. of 49th {IEEE} Conf. on Decision and Control}, Dec. 2010.
\newblock To appear.

\bibitem{kn:Gyorgy10}
G.~D{\'a}n and H.~Sandberg, ``Stealth attacks and protection schemes for state
  estimators in power systems,'' in {\em Proc. of {IEEE} SmartGridComm}, Oct.
  2010.

\bibitem{kn:Wu90_Survey}
F.~F. Wu, ``Power system state estimation: a survey,'' {\em Int. J. Elec. Power
  and Energy Systems}, Apr. 1990.

\bibitem{kn:MATPOWER}
R.~D. Zimmerman, C.~E. {Murillo-S\'{a}nchez}, and R.~J. Thomas, ``{MATPOWER's}
  extensible optimal power flow architecture,'' in {\em Power and Energy
  Society General Meeting}, pp.~1--7, IEEE, July 2009.

\end{thebibliography}

\end{document}